\pdfoutput=1
\documentclass[english,linesnumbered,ruled,vlined,11pt]{article}
\usepackage[a4paper]{geometry}
\usepackage{babel}
\usepackage{verbatim}
\usepackage{mathtools}
\usepackage{graphicx}
\usepackage{booktabs} 
\usepackage{enumitem}
\usepackage{bm}
\usepackage{algorithm, algpseudocode, multirow}
\usepackage{amsmath}
\usepackage{amsthm}
\usepackage{amssymb}
\usepackage{minitoc}
\usepackage{color,xcolor}
\usepackage{framed}
\usepackage[most]{tcolorbox}
\usepackage{listings}
\usepackage{xcolor}
\usepackage{wrapfig}

\lstdefinelanguage{julia}{
  morekeywords={abstract,break,case,catch,const,continue,do,else,elseif,end,export,false,for,function,immutable,import,importall,if,in,macro,module,quote,return,struct,true,try,type,typealias,using,while},
  sensitive=true,
  morecomment=[l]{\#},
  morestring=[b]",
}

\newtcblisting{codeblock}[2][]{ 
    arc=3pt,
    colback=gray!5,
    colframe=gray!30,
    listing only,
    title=#2, 
    fonttitle=\bfseries\small,
    coltitle=black,
    colbacktitle=gray!20,
    listing options={
        language=#2, 
        basicstyle=\ttfamily\small,
        keywordstyle=\color{blue}\bfseries,
        commentstyle=\color{green!50!black}\itshape,
        stringstyle=\color{purple},
        showstringspaces=false,
        breaklines=true,
        #1
    }, 
    enhanced,
    attach boxed title to top left={yshift=-2mm, xshift=2mm},
}
\definecolor{shadecolor}{rgb}{0.9,0.9,0.9} 
\usepackage{pifont} 
\usepackage[numbers]{natbib}
\usepackage[pdfusetitle,bookmarks=true,bookmarksnumbered=false,bookmarksopen=false,  breaklinks=false,pdfborder={0 0 1},backref=false,colorlinks=false]{hyperref}

\usepackage{cleveref}

\makeatletter
\providecommand{\theHALG@line}{}
\renewcommand{\theHALG@line}{\thealgorithm.\arabic{ALG@line}}

\theoremstyle{definition}
\newtheorem{definition}{\protect\definitionname}[section]
\theoremstyle{plain}

\theoremstyle{plain}

\theoremstyle{plain}

\theoremstyle{plain}

\theoremstyle{remark}

\theoremstyle{definition}

\providecommand{\examplename}{Example}
\theoremstyle{corollary}

\newcommand{\qbf}{\mathbf{q}}

\providecommand{\assumptionname}{Assumption}
\providecommand{\definitionname}{Definition}
\providecommand{\lemmaname}{Lemma}
\providecommand{\propositionname}{Proposition}
\providecommand{\remarkname}{Remark}
\providecommand{\theoremname}{Theorem}
\providecommand{\corollaryname}{Corollary}

\crefdefaultlabelformat{#2\textbf{#1}#3}
\crefname{section}{\textbf{section}}{\textbf{sections}}
\Crefname{section}{\textbf{Section}}{\textbf{Sections}}
\crefname{thm}{\textbf{Theorem}}{\textbf{theorems}}
\Crefname{thm}{\textbf{Theorem}}{\textbf{Theorems}}
\crefname{lem}{\textbf{Lemma}}{\textbf{lemmas}}
\Crefname{lem}{\textbf{Lemma}}{\textbf{Lemmas}}
\crefname{prop}{\textbf{proposition}}{\textbf{propositions}}
\Crefname{prop}{\textbf{Proposition}}{\textbf{Propositions}}
\crefname{algorithm}{\textbf{algorithm}}{\textbf{algorithms}}
\Crefname{algorithm}{\textbf{Algorithm}}{\textbf{Algorithms}}
\crefname{coro}{\textbf{Corollary}}{\textbf{corollaries}}
\Crefname{coro}{\textbf{Corollary}}{\textbf{corollaries}}
\crefname{defn}{\textbf{Definition}}{\textbf{definitions}}
\Crefname{defn}{\textbf{Definition}}{\textbf{definitions}}
\crefname{table}{\textbf{Table}}{\textbf{tables}}
\Crefname{table}{\textbf{Table}}{\textbf{tables}}
\crefname{figure}{\textbf{Figure}}{\textbf{figures}}
\Crefname{figure}{\textbf{Figure}}{\textbf{figures}}
\crefname{exple}{\textbf{Example}}{\textbf{examples}}
\Crefname{exple}{\textbf{Example}}{\textbf{examples}}
\Crefname{assumption}{\textbf{Assumption}}{\textbf{Assumptions}}
\crefname{assumption}{\textbf{Assumption}}{\textbf{Assumptions}}
\Crefname{rem}{\textbf{Remark}}{\textbf{Remarks}}
\crefname{rem}{\textbf{Remark}}{\textbf{Remarks}}

\usepackage{fancybox,ulem,subfig}
\tcbuselibrary{breakable,theorems,skins}

\newtcbtheorem[number within=section]{exbox}{Example}{
  enhanced,
  colback=gray!3,           
  colframe=gray!70!black,   
  boxrule=0.6pt,
  arc=2mm,                  
  left=8pt,right=8pt,top=8pt,bottom=8pt,
  fonttitle=\bfseries,       
  attach boxed title to top left={yshift=-1mm, xshift=2mm},
}{ex}

\makeatother

\begin{document}
\global\long\def\inprod#1#2{\left\langle #1,#2\right\rangle }%

\global\long\def\inner#1#2{\left\langle#1,#2\right\rangle}%

\global\long\def\binner#1#2{\big\langle#1,#2\big\rangle}%

\global\long\def\norm#1{\left\|#1\right\|}%

\global\long\def\bnorm#1{\big\Vert#1\big\Vert}%

\global\long\def\Bnorm#1{\Big\Vert#1\Big\Vert}%

\global\long\def\brbra#1{\left(#1\right)}%

\global\long\def\Brbra#1{\left(#1\right)}%

\global\long\def\rbra#1{(#1)}%

\global\long\def\sbra#1{[#1]}%

\global\long\def\bsbra#1{\big[#1\big]}%

\global\long\def\Bsbra#1{\Big[#1\Big]}%

\global\long\def\abs#1{\left|#1\right|}%

\global\long\def\babs#1{\big\vert#1\big\vert}%

\global\long\def\cbra#1{\{#1\}}%

\global\long\def\bcbra#1{\left\{#1\right\}}%

\global\long\def\Bcbra#1{\Big\{#1\Big\}}%

\global\long\def\vertiii#1{\left\vert \kern-0.25ex  \left\vert \kern-0.25ex  \left\vert #1\right\vert \kern-0.25ex  \right\vert \kern-0.25ex  \right\vert }%

\global\long\def\matr#1{\bm{#1}}%

\global\long\def\til#1{\tilde{#1}}%

\global\long\def\wtil#1{\widetilde{#1}}%

\global\long\def\wh#1{\widehat{#1}}%

\global\long\def\mcal#1{\mathcal{#1}}%

\global\long\def\mbb#1{\mathbb{#1}}%

\global\long\def\mtt#1{\mathtt{#1}}%

\global\long\def\ttt#1{\texttt{#1}}%

\global\long\def\dtxt{\textrm{d}}%

\global\long\def\aeq{\overset{(a)}{=}}%

\global\long\def\bignorm#1{\bigl\Vert#1\bigr\Vert}%

\global\long\def\Bignorm#1{\Bigl\Vert#1\Bigr\Vert}%

\global\long\def\rmn#1#2{\mathbb{R}^{#1\times#2}}%

\global\long\def\deri#1#2{\frac{d#1}{d#2}}%

\global\long\def\pderi#1#2{\frac{\partial#1}{\partial#2}}%

\global\long\def\limk{\lim_{k\rightarrow\infty}}%

\global\long\def\trans{\textrm{T}}%

\global\long\def\onebf{\mathbf{1}}%

\global\long\def\zerobf{\mathbf{0}}%

\global\long\def\zero{\bm{0}}%

\global\long\def\Euc{\mathrm{E}}%

\global\long\def\Expe{\mathbb{E}}%

\global\long\def\rank{\mathrm{rank}}%

\global\long\def\range{\mathrm{range}}%

\global\long\def\diam{\mathrm{diam}}%

\global\long\def\epi{\mathrm{epi} }%

\global\long\def\inte{\operatornamewithlimits{int}}%

\global\long\def\dist{\operatornamewithlimits{dist}}%

\global\long\def\proj{\operatorname{Proj}}%

\global\long\def\cov{\mathrm{Cov}}%

\global\long\def\argmin{\operatorname{argmin}}%

\global\long\def\argmax{\operatornamewithlimits{argmax}}%

\global\long\def\where{\operatornamewithlimits{where}}%

\global\long\def\tr{\operatornamewithlimits{tr}}%

\global\long\def\dis{\operatornamewithlimits{dist}}%

\global\long\def\sign{\operatornamewithlimits{sign}}%

\global\long\def\prob{\mathrm{Prob}}%

\global\long\def\st{\operatornamewithlimits{s.t.}}%

\global\long\def\dom{\mathrm{dom}}%

\global\long\def\prox{\mathrm{prox}}%

\global\long\def\diag{\mathrm{diag}}%

\global\long\def\and{\mathrm{and}}%

\global\long\def\as{\textup{a.s.}}%

\global\long\def\ae{\textup{a.e.}}%

\global\long\def\Var{\operatornamewithlimits{Var}}%

\global\long\def\Cov{\operatornamewithlimits{Cov}}%

\global\long\def\raw{\rightarrow}%

\global\long\def\law{\leftarrow}%

\global\long\def\Raw{\Rightarrow}%

\global\long\def\Law{\Leftarrow}%

\global\long\def\vep{\varepsilon}%

\global\long\def\eps{\varepsilon}%

\global\long\def\dom{\operatornamewithlimits{dom}}%

\global\long\def\err{\text{err}}%

\global\long\def\soc{\operatorname{soc}}%

\global\long\def\rsoc{\operatorname{rsoc}}%

\global\long\def\tsum{{\textstyle {\sum}}}%

\global\long\def\Cbb{\mathbb{C}}%

\global\long\def\Ebb{\mathbb{E}}%

\global\long\def\Fbb{\mathbb{F}}%

\global\long\def\Nbb{\mathbb{N}}%

\global\long\def\Rbb{\mathbb{R}}%

\global\long\def\extR{\widebar{\mathbb{R}}}%

\global\long\def\Pbb{\mathbb{P}}%

\global\long\def\Mrm{\mathrm{M}}%

\global\long\def\Acal{\mathcal{A}}%

\global\long\def\Bcal{\mathcal{B}}%

\global\long\def\Ccal{\mathcal{C}}%

\global\long\def\Dcal{\mathcal{D}}%

\global\long\def\Ecal{\mathcal{E}}%

\global\long\def\Fcal{\mathcal{F}}%

\global\long\def\Gcal{\mathcal{G}}%

\global\long\def\Hcal{\mathcal{H}}%

\global\long\def\Ical{\mathcal{I}}%

\global\long\def\Kcal{\mathcal{K}}%
\global\long\def\calK{\mathcal{K}}%
\global\long\def\rtcomment#1{\hfill\textit{#1}}%

\global\long\def\Lcal{\mathcal{L}}%

\global\long\def\Mcal{\mathcal{M}}%

\global\long\def\Ncal{\mathcal{N}}%

\global\long\def\Ocal{\mathcal{O}}%

\global\long\def\Pcal{\mathcal{P}}%

\global\long\def\Scal{\mathcal{S}}%

\global\long\def\Tcal{\mathcal{T}}%

\global\long\def\Xcal{\mathcal{X}}%

\global\long\def\Ycal{\mathcal{Y}}%

\global\long\def\Zcal{\mathcal{Z}}%

\global\long\def\i{i}%

\global\long\def\abf{\mathbf{a}}%

\global\long\def\Nbf{\mathbf{N}}%

\global\long\def\bbf{\mathbf{b}}%

\global\long\def\cbf{\mathbf{c}}%

\global\long\def\fbf{\mathbf{f}}%

\global\long\def\sbf{\mathbf{s}}%

\global\long\def\qbf{\mathbf{q}}%

\global\long\def\gbf{\mathbf{g}}%

\global\long\def\lambf{\bm{\lambda}}%

\global\long\def\alphabf{\bm{\alpha}}%

\global\long\def\sigmabf{\bm{\sigma}}%

\global\long\def\thetabf{\bm{\theta}}%

\global\long\def\deltabf{\bm{\delta}}%

\global\long\def\lbf{\mathbf{l}}%

\global\long\def\ubf{\mathbf{u}}%

\global\long\def\pbf{\mathbf{\mathbf{p}}}%

\global\long\def\vbf{\mathbf{v}}%

\global\long\def\rbf{\mathbf{r}}%

\global\long\def\wbf{\mathbf{w}}%

\global\long\def\xbf{\mathbf{x}}%

\global\long\def\tbf{\mathbf{t}}%

\global\long\def\ybf{\mathbf{y}}%

\global\long\def\zbf{\mathbf{z}}%

\global\long\def\dbf{\mathbf{d}}%

\global\long\def\hbf{\mathbf{h}}%

\global\long\def\Wbf{\mathbf{W}}%

\global\long\def\Abf{\mathbf{A}}%

\global\long\def\Ubf{\mathbf{U}}%

\global\long\def\Pbf{\mathbf{P}}%

\global\long\def\Gbf{\mathbf{G}}%

\global\long\def\Ibf{\mathbf{I}}%

\global\long\def\Ebf{\mathbf{E}}%

\global\long\def\ebf{\mathbf{e}}%

\global\long\def\Mbf{\mathbf{M}}%

\global\long\def\Dbf{\mathbf{D}}%

\global\long\def\Qbf{\mathbf{Q}}%

\global\long\def\Lbf{\mathbf{L}}%

\global\long\def\Pbf{\mathbf{P}}%

\global\long\def\Xbf{\mathbf{X}}%

\global\long\def\abm{\bm{a}}%

\global\long\def\bbm{\bm{b}}%

\global\long\def\cbm{\bm{c}}%

\global\long\def\dbm{\bm{d}}%

\global\long\def\ebm{\bm{e}}%

\global\long\def\fbm{\bm{f}}%

\global\long\def\gbm{\bm{g}}%

\global\long\def\hbm{\bm{h}}%

\global\long\def\pbm{\bm{p}}%

\global\long\def\qbm{\bm{q}}%

\global\long\def\rbm{\bm{r}}%

\global\long\def\sbm{\bm{s}}%

\global\long\def\tbm{\bm{t}}%

\global\long\def\ubm{\bm{u}}%

\global\long\def\vbm{\bm{v}}%

\global\long\def\wbm{\bm{w}}%

\global\long\def\xbm{\bm{x}}%

\global\long\def\ybm{\bm{y}}%

\global\long\def\zbm{\bm{z}}%

\global\long\def\Abm{\bm{A}}%

\global\long\def\Bbm{\bm{B}}%

\global\long\def\Cbm{\bm{C}}%

\global\long\def\Dbm{\bm{D}}%

\global\long\def\Ebm{\bm{E}}%

\global\long\def\Fbm{\bm{F}}%

\global\long\def\Gbm{\bm{G}}%

\global\long\def\Hbm{\bm{H}}%

\global\long\def\Ibm{\bm{I}}%

\global\long\def\Jbm{\bm{J}}%

\global\long\def\Lbm{\bm{L}}%

\global\long\def\Obm{\bm{O}}%

\global\long\def\Pbm{\bm{P}}%

\global\long\def\Qbm{\bm{Q}}%

\global\long\def\Rbm{\bm{R}}%

\global\long\def\Ubm{\bm{U}}%

\global\long\def\Vbm{\bm{V}}%

\global\long\def\Wbm{\bm{W}}%

\global\long\def\Xbm{\bm{X}}%

\global\long\def\Ybm{\bm{Y}}%

\global\long\def\Zbm{\bm{Z}}%

\global\long\def\lambm{\bm{\lambda}}%

\global\long\def\alphabm{\bm{\alpha}}%

\global\long\def\albm{\bm{\alpha}}%

\global\long\def\taubm{\bm{\tau}}%

\global\long\def\mubm{\bm{\mu}}%

\global\long\def\yrm{\mathrm{y}}%

\global\long\def\rone{\text{\textrm{I}}}%

\global\long\def\rtwo{\text{II}}%

\global\long\def\rthree{\text{\textrm{III}}}%

\global\long\def\rfour{\text{\textrm{IV}}}%

\global\long\def\rfive{\text{V}}%

\global\long\def\rsix{\text{\textrm{VI}}}%

\global\long\def\rseven{\text{VI\textrm{I}}}%

\global\long\def\reight{\text{VI\textrm{I}I}}%

\global\long\def\PDCS{\text{PDCS}}%

\global\long\def\vbfp{\vbf_{\text{p}}}%

\global\long\def\vbfd{\vbf_{\text{d}}}%

\global\long\def\tp{t_{\text{p}}}%

\global\long\def\td{t_{\text{d}}}%

\global\long\def\coneFamily{\text{Projection-friendly Cone}}

\global\long\def\vpr{v_{\textrm{pr}}}%
\global\long\def\vps{v_{\textrm{ps}}}%
\global\long\def\vpt{v_{\textrm{pt}}}%
\global\long\def\vdr{v_{\textrm{dr}}}%
\global\long\def\vds{v_{\textrm{ds}}}%
\global\long\def\vdt{v_{\textrm{dt}}}%
\global\long\def\dr{d_{\textrm{r}}}%
\global\long\def\ds{d_{\textrm{s}}}%
\global\long\def\dt{d_{\textrm{t}}}%

\global\long\def\bgeq{\overset{(b)}{\geq}}%
\global\long\def\sp{s_{\textrm{p}}}%
\global\long\def\rd{r_{\textrm{d}}}%
\global\long\def\tildesp{\tilde{s}_{\textrm{p}}}%
\global\long\def\tilderd{\tilde{r}_{\textrm{d}}}%

\newcommand{\qedtext}{{\unskip\nobreak\hfil\penalty50\hskip1em\hbox{}\nobreak\hfil$\square$\parfillskip=0pt\finalhyphendemerits=0\par}}
 
\providecommand{\tabularnewline}{\\}
\floatstyle{ruled}
\newfloat{algorithm}{tbp}{loa}
\providecommand{\algorithmname}{Algorithm}
\floatname{algorithm}{\protect\algorithmname}

\title{A Technical Note on the Implementation and Use of PDCS}

\author{Zhenwei Lin\thanks{lin2193@purdue.edu, School of Industrial Engineering, Purdue University} \qquad Zikai Xiong\thanks{zxiong84@gatech.edu, H. Milton Stewart School of Industrial and Systems Engineering, Georgia Institute of Technology} \qquad Dongdong Ge\thanks{ddge@sjtu.edu.cn, Antai School of Economics and Management, Shanghai Jiao Tong University} \qquad Yinyu Ye\thanks{yyye@stanford.edu, Department of Management Science and Engineering and Institute of Computational \& Mathematical Engineering, Stanford University}  }

\maketitle
\begin{abstract}
  This technical note documents the implementation and use of the Primal-Dual Conic Programming Solver (PDCS), a first-order solver for large-scale conic optimization problems introduced by Lin et al. (arXiv:2505.00311). It describes the algorithmic and implementation details underlying PDCS, including the restarted primal-dual hybrid gradient method framework, adaptive step-size selection, adaptive reflected Halpern iterations, adaptive restarts, and diagonal preconditioning. It also provides practical instructions for using PDCS, including its interfaces with JuMP and CVXPY, solver options, and illustrative code examples. PDCS is available at \url{https://github.com/ZikaiXiong/PDCS}  under the Apache License 2.0.   
\end{abstract}
\section{Introduction}
Convex conic programs arise in finance, machine learning, energy systems, and control. While interior-point methods are robust for moderate-sized problems, their reliance on matrix factorizations limits their scalability for large-scale instances involving millions of variables or constraints. In contrast, first-order methods rely on matrix-vector multiplications, which are memory-efficient and well suited to parallel computation on graphics processing units (GPUs).

PDCS was introduced by \citet{lin2025practical}, who prove the convergence of its base algorithm, characterize the projections onto diagonally rescaled cones, and compare cuPDCS with other solvers. This note documents the implementation of PDCS, a matrix-free solver for conic optimization problems whose cones are Cartesian products of basic cone blocks. Built on the PDHG framework, PDCS incorporates several practical enhancements, including adaptive restarting based on the normalized duality gap, adaptive reflected Halpern iteration, primal-weight updates, and diagonal preconditioning. The note also describes the GPU implementation and the available user interfaces.

\subsection{\label{sec:preliminaries}Preliminaries on Conic Optimization and PDHG} 

We consider the following conic program:
\begin{equation}\label{eq:prob_cp}
\small
    \min_{\substack{\xbf=(\xbf_{1},\xbf_{2}):  \xbf_{1}\in\mbb R^{n_{1}},\xbf_{2}\in\mbb R^{n_{2}}}}\inner{\cbf}{\xbf}\ \ \st\ \mathbf{G}\xbf-\mathbf{h}\in\mcal K_{\textrm{d}}^{*}\ ,\ \mathbf{l}\leq\xbf_{1}\leq\mathbf{u}\ ,\ \xbf_{2}\in\mcal K_{\textrm{p}}\ ,
\end{equation} 
in which,
\begin{itemize}
    \item $\mathbf{G}\in \mbb R^{m\times n}$ is the constraint matrix, 
    \item $\mathbf{l}\in (\mathbb R\cup\{-\infty\})^{n_1}$ and $\mathbf{u}\in (\mathbb R\cup\{\infty\})^{n_1}$ denote the lower and upper bounds on the components of $\mathbf{x}_1$,
    \item The set $\mathcal K_{\textrm p}$ denotes the primal cone, which is a Cartesian product of cone blocks, say $\mathcal K_{\textrm p,1}\times\cdots\times\mathcal K_{\textrm p,l_{\textrm p}}$, where $l_{\textrm p}$ is the number of cone blocks. In PDCS, these blocks may include the zero cone $\zerobf^d$, the nonnegative orthant $\mbb R_+^d$, the second-order cone $\mcal K_{\soc}^{d+1}$, the exponential cone $\mcal K_{\exp}$, and the dual exponential cone $\mcal K_{\exp}^*$,
    \item $\mcal K_{\textrm{d}}^{*}$ denotes the dual cone of $\mcal K_{\textrm{d}}$.
\end{itemize}

The dual problem of \eqref{eq:prob_cp} is given by:
\begin{equation}\label{eq:dual_conic}
    \begin{aligned}
\max_{\substack{\ybf,\  \bm{\lambda} = (\bm{\lambda}_{1},\bm{\lambda}_{2}):
\\
\ybf\in \mbb R^m, \ \bm{\lambda}_{1}\in\mathbb{R}^{n_1}, \ \bm{\lambda}_{2}\in\mathbb{R}^{n_2}}}  \ \langle \ybf, \mathbf{h}\rangle + \lbf^{\top}\bm{\lambda}_{1}^{+}-\ubf^{\top}\bm{\lambda}_{1}^{-}\ \ \st  \cbf-\Gbf^{\top}\ybf=\bm{\lambda}\ ,\ \ybf\in\mcal K_{\textrm{d}}, \ \bm{\lambda}_{1}\in\Lambda\ ,\ \bm{\lambda}_{2}\in\mcal K_{\textrm{p}}^{*}\ ,
    \end{aligned}
\end{equation} 
where \begin{equation}\label{eq:defn_Lambda1}
    \Lambda := \Lambda_1\times \Lambda_2 \times \cdots \times \Lambda_{n_1} 
    \ \ \text{ and } \ \
    \Lambda_i :=\begin{cases}
        \{0\}\ , & \text{if }\lbf_{i}=-\infty,\ubf_{i}=+\infty \\
        \mbb R_{-}\ , & \text{if }\lbf_{i}=-\infty,\ubf_{i}\in\mbb R \\
        \mbb R_{+}\ , & \text{if }\lbf_{i}\in\mbb R,\ubf_{i}=+\infty \\
        \mbb R\ , & \text{otherwise}
        \end{cases} \text{ for }i=1,2,\dots,n_1 \ .
\end{equation}
Similarly, $\mcal K_{\textrm d}=\mcal K_{\textrm d,1}\times\cdots\times\mcal K_{\textrm d,l_{\textrm d}}$ is a Cartesian product of cone blocks, and its dual cone is $\mcal K_{\textrm d}^*=\mcal K_{\textrm d,1}^*\times\cdots\times\mcal K_{\textrm d,l_{\textrm d}}^*$.

The problem~\eqref{eq:prob_cp} has an equivalent primal-dual formulation, expressed as the following saddle-point problem on the Lagrangian $\mcal L(\xbf,\ybf)$: 
\begin{equation}\label{eq:pd_prob}
    \min_{\xbf \in \mcal X}  \ \max_{\ybf \in \mcal Y}\ \ \mcal L(\xbf,\ybf):=\inner{\cbf}{\xbf} -\inner{\ybf}{\mathbf{G}\xbf-\mathbf{h}}\ ,
\end{equation} 
where $\mcal X:=[\lbf,\ubf]\times \mcal K_{\textrm{p}}$ and $\mcal Y:=\mcal K_{\textrm{d}}$. A saddle point $(\xbf,\ybf)$ of~\eqref{eq:pd_prob} corresponds to an optimal primal solution $\xbf$ for \eqref{eq:prob_cp} and an optimal dual solution $(\ybf, \cbf-\Gbf^\top \ybf)$  for \eqref{eq:dual_conic}.  
\paragraph{Optimality measures.}
For $\bm{\lambda}=\cbf-\Gbf^{\top}\ybf$, the primal infeasibility, dual infeasibility, and objective gap are measured by
\begin{equation}\label{eq:tolerance-3}
\begin{aligned}
\err_{\textrm{abs,p}}&:=\|\Gbf\xbf-\hbf-\proj_{\mcal K_{\textrm d}^{*}}(\Gbf\xbf-\hbf)\|_2,\\
\err_{\textrm{abs,d}}&:=\|\bm{\lambda}-\proj_{\Lambda\times\mcal K_{\textrm p}^{*}}(\bm{\lambda})\|_2,\\
\err_{\textrm{abs,gap}}&:=\bigl|\cbf^\top\xbf-\ybf^\top\hbf-\lbf^\top\bm{\lambda}_1^++\ubf^\top\bm{\lambda}_1^-\bigr|.
\end{aligned}
\end{equation}

\paragraph{Stopping criterion.}
For termination, the GPU solver uses the same formulas with every vector $2$-norm replaced by the infinity norm. We denote these residual measures by a superscript $(\infty)$, while the objective-gap measures remain unchanged. Specifically,
\begin{equation}\label{eq:stopping_errors}
\begin{aligned}
\err_{\textrm{rel,p}}^{(\infty)}&:=\frac{\err_{\textrm{abs,p}}^{(\infty)}}{1+\max\{\|\Gbf\xbf\|_\infty,\|\hbf\|_\infty,\|\proj_{\mcal K_{\textrm d}^{*}}(\Gbf\xbf-\hbf)\|_\infty\}},\\
\err_{\textrm{rel,d}}^{(\infty)}&:=\frac{\err_{\textrm{abs,d}}^{(\infty)}}{1+\max\{\|\cbf\|_\infty,\|\Gbf^\top\ybf\|_\infty\}},\\
\err_{\textrm{rel,gap}}&:=\frac{\err_{\textrm{abs,gap}}}{1+\max\{\lvert\cbf^\top\xbf\rvert,\lvert\ybf^\top\hbf+\lbf^\top\bm{\lambda}_1^+-\ubf^\top\bm{\lambda}_1^-\rvert\}}.
\end{aligned}
\end{equation}
The solver declares an optimal solution when either
\begin{equation}\label{eq:stopping}
\begin{aligned}
&\max\{\err_{\textrm{abs,p}}^{(\infty)},\err_{\textrm{abs,d}}^{(\infty)},\err_{\textrm{abs,gap}}\}<\varepsilon_{\textrm{abs}},\\
&\text{or}\quad\max\{\err_{\textrm{rel,p}}^{(\infty)},\err_{\textrm{rel,d}}^{(\infty)},\err_{\textrm{rel,gap}}\}<\varepsilon_{\textrm{rel}},
\end{aligned}
\end{equation}
where $\varepsilon_{\textrm{abs}}$ and $\varepsilon_{\textrm{rel}}$ are set by \texttt{abs\_tol} and \texttt{rel\_tol}. All errors are evaluated in the original, unscaled problem. 
    
\subsection{\label{subsec:PDHG_CP}PDHG for conic programs}
PDCS uses PDHG as its base algorithm. Let $\tau$ and $\sigma$ denote the primal and dual step sizes, respectively. One iteration of PDHG for solving~\eqref{eq:pd_prob} from iterate $\zbf = (\xbf,\ybf)$, denoted by $\text{OnePDHG}(\zbf)$, is given by:
\begin{equation}
    \hat{\zbf}=\texttt{OnePDHG}(\zbf):=\begin{cases}
    \ensuremath{\hat{\xbf}=\proj_{[\lbf,\ubf]\times\mcal K_{\textrm{p}}}\bcbra{\xbf-\tau(\cbf-\Gbf^{\top}\ybf)}}\\
    \ensuremath{\hat{\ybf}=\proj_{\mcal K_{\textrm{d}}}\bcbra{\ybf+\sigma(\hbf-\Gbf(2\hat{\xbf}-\xbf))}}
    \end{cases}.\label{eq:PDHG}
\end{equation} 
When $\tau\sigma\|\Gbf\|^2<1$, the iterates generated by~\eqref{eq:PDHG} converge to a saddle point of~\eqref{eq:pd_prob}. See, for example, \citet{chambolle2011first,lu2022infimal,xiong2024role}. 

\section{\label{sec:Practical-enhancement-technique}
Algorithmic enhancements in PDCS}

PDCS combines adaptive step sizes (Section~\ref{sec:adaptive_stepsize}), reflected Halpern iteration (Section~\ref{sec:adaptive_reflection}), adaptive restarts (Section~\ref{sec:adaptive_restart}), and primal weight updates (Section~\ref{sec:weight_adjustment}), as summarized in Algorithm~\ref{alg:PDCS}. Interface-specific defaults are described in Section~\ref{sec:usage}.
\begin{algorithm}[htbp]
\caption{Primal-Dual Conic Programming Solver (PDCS) without preconditioning \label{alg:PDCS}}
    \begin{algorithmic}[1]
        \Require{Initial iterate $\bar{\zbf}^{0,0}=\zbf^{0,0}=(\xbf^{0,0},\ybf^{0,0})$, global trial counter $\bar{k}=1$, initial step size $\hat{\eta}^{0,0}=1/\norm{\Gbf}_{\infty}$, $t\leftarrow 0$, initial primal weight $\omega^0\leftarrow \texttt{InitializePrimalWeight}(\cbf,\hbf) \quad$ (Section~\ref{sec:weight_adjustment}).}
    \Repeat
    \State{$k\leftarrow 0$, $\zbf_{c}^{t,0}\leftarrow\zbf^{t,0}$}
    \While{neither the restart condition (Section~\ref{sec:adaptive_restart}) nor the termination criterion holds}
        \State{$\hat{\zbf}^{t,k+1}, \eta^{t,k+1}, \hat{\eta}^{t,k+1}, \bar{k}\leftarrow$ \texttt{AdaptiveStepPDHG}($\zbf^{t,k}, \omega^t, \hat{\eta}^{t,k},\bar{k}$)\label{alg_pdcs:line4}\rtcomment{Section~\ref{sec:adaptive_stepsize}} }
        \State{$\alpha^{t,k}\leftarrow$ \texttt{AdaptiveReflectionParameter}($\hat{\zbf}^{t,k+1}$)\label{alg_pdcs:line5}\rtcomment{Section~\ref{sec:adaptive_reflection}}}
        \State{$\zbf^{t,k+1}\leftarrow \tfrac{k+1}{k+2}\big((1+\alpha^{t,k})\hat{\zbf}^{t,k+1}-\alpha^{t,k}\zbf^{t,k}\big) + \frac{1}{k+2}\zbf^{t,0}$\label{alg:line:ref_Halpern}}
        \State{$\bar{\zbf}^{t,k+1}\leftarrow \tsum_{i=1}^{k+1}\eta^{t,i}\zbf^{t,i}/{\tsum_{i=1}^{k+1}\eta^{t,i}}$\label{alg:line:average}} 
        \State{$\zbf_{c}^{t,k+1}\leftarrow\texttt{GetRestartCandidate}(\zbf^{t,k+1},\bar{\zbf}^{t,k+1})$\label{alg:line:GetRestartCandidate}\rtcomment{Section~\ref{sec:adaptive_restart}}}
        \State{$k\leftarrow k+1$}
    \EndWhile
    \State{$\zbf^{t+1,0}\leftarrow \zbf_c^{t,k}$, $\omega^{t+1}\leftarrow \texttt{PrimalWeightUpdate}(\zbf^{t+1,0}, \zbf^{t,0},\omega^{t})$\rtcomment{Section~\ref{sec:weight_adjustment}}}
    \State{$\hat{\eta}^{t+1,0}\leftarrow\hat{\eta}^{t,k}$, $t\leftarrow t+1$}
    \Until{termination criteria hold}
    \State{\textbf{Return} $\zbf^{t,0}$.}
    \end{algorithmic}
\end{algorithm}  
In addition, PDCS also applies diagonal preconditioning to improve the conditioning of problem instances. This preconditioning procedure is described in Section~\ref{subsec:Projection-onto-convex}. 

In Algorithm~\ref{alg:PDCS}, the function \texttt{AdaptiveStepPDHG} in Line~\ref{alg_pdcs:line4} performs a line search to determine an appropriate step size, and the function \texttt{AdaptiveReflectionParameter} in Line~\ref{alg_pdcs:line5} computes the reflection coefficient. Lines~\ref{alg:line:ref_Halpern} and~\ref{alg:line:average} perform a reflected Halpern iteration and weighted averaging of iterates, and the function \texttt{GetRestartCandidate} in Line~\ref{alg:line:GetRestartCandidate} selects a candidate for restarting based on either the normalized duality gap or KKT error. At each outer iteration, the primal weight $\omega$ is updated based on the progress of primal and dual iterates, and this weight balances the primal and dual step sizes in \texttt{AdaptiveStepPDHG}. 

\subsection{\label{sec:adaptive_stepsize}Adaptive step size} 

PDHG converges when $\tau\sigma\|\Gbf\|^2<1$, where $\|\Gbf\|$ is the spectral norm~\citep{chambolle2011first}. Because this bound can be conservative in practice, PDCS uses a line-search heuristic similar to that of PDLP~\citep{applegate2021practical}:
\begin{equation}\label{eq:adaptive_step_size}
    \tau = \frac{\eta}{\omega},\quad \sigma = \eta \omega,\quad \eta \leq \frac{\norm{\zbf^{t+1}-\zbf^t}_{\omega}^2}{2\abs{(\ybf^{t+1}-\ybf^t)^\top \Gbf (\xbf^{t+1}-\xbf^t)}},
\end{equation}
where $\norm{\zbf}_{\omega} := \sqrt{\omega\norm{\xbf}^2 + \frac{\norm{\ybf}^2}{\omega}}$, and $\omega$ is the primal weight. 

The update of $\eta$ is carried out by the procedure described in Function~\ref{alg:adaptive_pdhg}. Here, $\bar{k}$ counts all PDHG steps, including rejected trials, starting from $\bar{k}=1$, and controls how aggressively $\eta$ is adjusted \citep{applegate2021practical}.

{\renewcommand{\algorithmname}{Function}\renewcommand{\thealgorithm}{1}\def\theHalgorithm{function.1}
\makeatletter\renewcommand{\theHALG@line}{function.\arabic{ALG@line}}\makeatother
\begin{algorithm}[h]
\caption{AdaptiveStepPDHG(($\xbf$,$\ybf$), $\omega$, $\eta$, $\bar{k}$)\label{alg:adaptive_pdhg}}
    \small
    \begin{algorithmic}[1]
    \While{True}
    \State{$\hat{\xbf} = \proj_{[\lbf,\ubf] \times \mcal K_{\textrm{p}}}\bcbra{\xbf - \frac{\eta}{\omega}(\cbf-\Gbf^\top \ybf)}$\label{alg_adaptive_pdhg_line2}}
    \State{$\hat{\ybf}=\proj_{\mcal K_{\textrm{d}}}\bcbra{\ybf + \eta \omega(\hbf-\Gbf(2\hat{\xbf}-\xbf))}$\label{alg_adaptive_pdhg_line3}}
    \State{$\bar{\eta}\leftarrow \frac{1}{2}{\norm{(\hat{\xbf}-\xbf,\hat{\ybf}-\ybf)}_{\omega}^2}/{\abs{(\hat{\ybf}-\ybf)^\top \Gbf (\hat{\xbf}-\xbf)}}$}
    \State{$\eta^{\prime}\leftarrow \min\bcbra{(1-(\bar{k}+1)^{-0.3})\bar{\eta}, (1+(\bar{k}+1)^{-0.6})\eta}$}
    \If{$\eta<\bar{\eta}$}
    \State{$\bar{k}\leftarrow\bar{k}+1$}
    \State{\Return{$(\hat{\xbf},\hat{\ybf}),\eta,\eta^{\prime},\bar{k}$}}
    \EndIf
    \State{$\eta\leftarrow \eta^{\prime},\bar{k}\leftarrow\bar{k}+1$}
    \EndWhile
    \end{algorithmic}
\end{algorithm}\addtocounter{algorithm}{-1}
}

\subsection{\label{sec:adaptive_reflection}Adaptive reflected Halpern iteration}  

Line~\ref{alg:line:ref_Halpern} in Algorithm~\ref{alg:PDCS} implements a reflected Halpern iteration. A similar technique has been used in practical first-order LP solvers such as HPR-LP~\citep{chen2024hpr} and r$^2$HPDHG~\citep{lu2024restarted}, where it accelerates convergence in practice. The standard Halpern iteration involves averaging the current iterate $\hat{\zbf}^{t,k+1}$ with the initial point $\zbf^{t,0}$ of the current inner loop. This anchoring mechanism mitigates oscillations in PDHG iterates and has been shown to yield improved worst-case convergence rates~\citep{lieder2021convergence}. In this setting, the point $\zbf^{t,0}$ serves as a fixed anchor. Originally introduced by \citet{halpern1967fixed}, Halpern iteration has recently been recognized for its role in accelerating a range of minimax optimization and fixed-point methods~\citep[see, e.g.,][]{yoon2025accelerated}.

The idea behind the reflected Halpern iteration is to apply this anchoring technique not directly to $\hat{\zbf}^{t,k+1}$, but to an extrapolated point between the newly computed iterate $\hat{\zbf}^{t,k+1}$ and the previous iterate ${\zbf}^{t,k}$. Specifically, Line~\ref{alg:line:ref_Halpern} replaces $\hat{\zbf}^{t,k+1}$ with $(1 + \alpha^{t,k})\hat{\zbf}^{t,k+1} - \alpha^{t,k}\zbf^{t,k}$, where $\alpha$ is the reflection parameter. For $\alpha=1$, this is the reflected operator $2\,\texttt{OnePDHG}-\Ibf$ applied to $\zbf^{t,k}$~\citep{ryu2016primer}. For linear programs, \citet{lu2024restarted} show that adding this reflection to restarted Halpern PDHG retains its theoretical guarantees and improves the complexity bound by a factor of two.

For reflection, the relative $2$-norm residuals are
\begin{equation}\label{eq:tolerance-2}
\begin{aligned}
\err_{\textrm{rel,p}}^{(2)}&:=\frac{\err_{\textrm{abs,p}}}{1+\max\{\|\Gbf\xbf\|_2,\|\hbf\|_2,\|\proj_{\mcal K_{\textrm d}^{*}}(\Gbf\xbf-\hbf)\|_2\}},\\
\err_{\textrm{rel,d}}^{(2)}&:=\frac{\err_{\textrm{abs,d}}}{1+\max\{\|\cbf\|_2,\|\Gbf^\top\ybf\|_2\}}.
\end{aligned}
\end{equation}
The reflection rule uses their maximum together with the relative objective gap in~\eqref{eq:stopping_errors},
\begin{equation}\label{eq:tolerance}
\text{maxErr}(\zbf):=\max\{\err_{\textrm{rel,p}}^{(2)},\err_{\textrm{rel,d}}^{(2)},\err_{\textrm{rel,gap}}\}.
\end{equation}

In practice, larger $\alpha$ yields more aggressive extrapolation and can accelerate convergence, but may also introduce instability and oscillations if set too high.
To balance acceleration and stability, PDCS adjusts the reflection coefficient adaptively based on the observed progress of the iterates. When the KKT error is large, indicating that the iterate is still far from optimality, a smaller coefficient is selected to promote conservative and stable updates. As the iterates begin to converge and the KKT error decreases, a larger coefficient is used to speed up progress. For conic programs, PDCS sets $\alpha^{t,k}$ from the KKT error at $\hat{\zbf}^{t,k+1}$ (Line~\ref{alg_pdcs:line5} of Algorithm~\ref{alg:PDCS}) by the rule
\begin{equation}\label{eq:adaptive_reflection}
    \texttt{AdaptiveReflectionParameter}(\zbf)=\proj_{[0,1]}\left(-0.1\cdot\log_{10}\left(\text{maxErr}(\zbf)\right)+0.2\right) \ ,
\end{equation}
where $\text{maxErr}$ is the relative Euclidean error measure defined in~\eqref{eq:tolerance}. This adjustment is separate from the infinity-norm stopping test in~\eqref{eq:stopping}.

\subsection{\label{sec:adaptive_restart}Adaptive restart} 
 
The restart mechanism resets the anchor point of the Halpern iteration to the restart candidate and updates the primal weight $\omega^t$  at the restart (see Section~\ref{sec:weight_adjustment}). It has been proven that PDHG with restarts (either average-iterate restarts or Halpern restarts) achieves linear convergence for LP problems~\citep{applegate2023faster,lu2024restarted}. For general conic programs, \citet{xiong2024role} establish the convergence of PDHG with average-iterate restarts, and \citet{lin2025practical} prove that the base algorithm of PDCS, reflected restarted Halpern PDHG with averaging, converges on general conic programs, linearly under a local error bound, when the step sizes, the reflection parameter, and the primal weight are fixed and restarts are decided by the fixed-point residual. The adaptive restart rule below is not covered by this analysis. Empirical studies have also shown that restart strategies can significantly improve the practical performance of PDHG for LP problems~\citep{applegate2021practical,applegate2023faster}.

To select the new anchor point, PDCS evaluates the normalized duality gap, which serves as a surrogate for optimality. This metric is defined as follows:
\begin{definition}[normalized duality gap~\citep{applegate2023faster}]\label{defn:normalized_dual_gap}
   For any $\zbf = (\xbf, \ybf) \in \mcal X \times \mcal Y$, $r > 0$, and step sizes $\tau > 0$, $\sigma > 0$ satisfying $\tau\sigma\|\Gbf\|^2<1$, define $\Mbf:=\bigg(\begin{smallmatrix}
\frac{1}{\tau}\Ibf & \Gbf^{\top}\\
\Gbf & \frac{1}{\sigma}\Ibf
\end{smallmatrix}\bigg)$ and $$\mcal{B} (r;\zbf):=\bcbra{\hat{\zbf}:=(\hat{\xbf},\hat{\ybf}):\hat{\xbf}\in \mcal X, \hat{\ybf}\in \mcal Y,\norm{\zbf-\hat{\zbf}}_{\Mbf}\leq r}\ ,$$ then the normalized duality gap is given by $
\rho(r;\zbf):=\frac{1}{r}\sup_{\hat{\zbf}\in \mcal B(r;\zbf)}[\mcal L(\xbf,\hat{\ybf})-\mcal L(\hat{\xbf},\ybf)].$
\end{definition}

The $\tau$ and $\sigma$ above correspond to the current step sizes $\tau^{t,k}$ and $\sigma^{t,k}$. 
\citet{applegate2023faster,xiong2024role} show that the above normalized duality gap $\rho(r;\zbf)$ serves as a measure of the optimality of $\zbf$. However, computing $\rho(r;\zbf)$  exactly involves a nontrivial projection onto $\mcal X \times \mcal Y$ under the $\Mbf$-norm, which is computationally expensive. To alleviate this, we approximate the duality gap by replacing $\Mbf$ with a diagonal surrogate matrix $\Nbf := \begin{pmatrix}
    \frac{1}{\tau}\Ibf &  \\
    & \frac{1}{\sigma}\Ibf
\end{pmatrix}$ and denote the resulting approximate gap as $\rho^{\Nbf}(\cdot,\cdot)$. This approximation $\rho^{\Nbf}(\cdot,\cdot)$ is equivalent to $\rho(\cdot,\cdot)$ up to a constant factor \citep{applegate2023faster,xiong2024role} and can be efficiently computed using a bisection search method \citep{applegate2023faster,xiong2024role}.
Details of this computation are given in Section~\ref{app:normalized_duality_gap}.

In PDCS, we compare the approximate normalized duality gaps of $\zbf^{t,k+1}$ and $\bar{\zbf}^{t,k+1}$, and select the one with the smaller gap as the new anchor point:
    \begin{equation*}
        \begin{aligned}
            & \texttt{GetRestartCandidate}(\zbf^{t,k+1},\bar{\zbf}^{t,k+1})
            \\
            & \quad \quad =\begin{cases}
        \zbf^{t,k+1} & \rho^{\Nbf}\brbra{\norm{\zbf^{t,k+1}-\zbf^{t,0}}_{\Nbf};\zbf^{t,k+1}}\leq\rho^{\Nbf}\brbra{\norm{\bar{\zbf}^{t,k+1}-\zbf^{t,0}}_{\Nbf};\bar{\zbf}^{t,k+1}}\\
        \bar{\zbf}^{t,k+1} & \text{otherwise}\ .
        \end{cases}
        \end{aligned}
        \end{equation*}
Let $\zbf_{c}^{t,k+1}$ denote the resulting restart candidate. A restart to $\zbf_{c}^{t,k+1}$ is triggered whenever its approximate normalized duality gap satisfies any of the following conditions. Here we use $\beta_{\text{sufficient}} = 0.4$, $\beta_{\text{necessary}} = 0.8$, and $\beta_{\text{artificial}} = 0.223$.  
\begin{enumerate}
    \item (Sufficient decay) $\rho^{\Nbf}\Brbra{\norm{\zbf_{c}^{t,k+1}-\zbf^{t,0}}_{\Nbf};\zbf_{c}^{t,k+1}}\leq\beta_\text{sufficient} \cdot \rho^{\Nbf}\brbra{\norm{\zbf^{t,0}-\zbf^{t-1,0}}_{\Nbf};\zbf^{t,0}}$
    \item (Necessary decay + no local progress)\newline $\rho^{\Nbf}\left(\norm{\zbf_{c}^{t,k+1}-\zbf^{t,0}}_{\Nbf};\zbf_{c}^{t,k+1}\right)>\rho^{\Nbf}\Brbra{\norm{\zbf_{c}^{t,k}-\zbf^{t,0}}_{\Nbf};\zbf_{c}^{t,k}}$ and\newline $\rho^{\Nbf}\Brbra{\norm{\zbf_{c}^{t,k+1} -\zbf^{t,0}}_{\Nbf};\zbf_{c}^{t,k+1}}\leq\beta_{\text{necessary}} \cdot \rho^{\Nbf}\brbra{\norm{\zbf^{t,0}-\zbf^{t-1,0}}_{\Nbf};\zbf^{t,0}}$
    \item (Long inner loop) $k \geq\beta_{\text{artificial}} \cdot \bar{k}$.
\end{enumerate}
Similar adaptive restart strategies have also been used in first-order LP solvers such as PDLP \citep{applegate2021practical}, cuPDLP \citep{lu2023cupdlp} and cuPDLP-c \citep{lu2023cupdlp-c}.

However, for general conic programs, numerical issues are more severe, particularly when evaluating $\rho^{\Nbf}$ for non-LP cones. In some cases, especially when the actual gap value is very small, the bisection method may produce negative outputs due to floating-point inaccuracies. When this occurs, PDCS switches to an alternative criterion based on the weighted KKT error: 
$$\text{KKT}_{\omega}(\zbf)=\sqrt{\omega^{2}\err_{\textrm{abs,p}}(\xbf)^{2}+\frac{1}{\omega^{2}}\err_{\text{\textrm{abs,d}}}(\bm{\lambda}_{1},\bm{\lambda}_{2})^{2}+\err_{\text{\textrm{abs,gap}}}(\xbf,\bm{\lambda}_{1},\bm{\lambda}_{2})^{2}}\ ,
$$
where $\err_{\textrm{abs,p}}, \err_{\textrm{abs,d}}, \err_{\textrm{abs,gap}}$ are as defined in~\eqref{eq:tolerance-3}, and the dual variables $\bm{\lambda}_1, \bm{\lambda}_2$ are recovered from $\ybf$ using the linear relation in~\eqref{eq:dual_conic}. When numerical issues happen, all instances of $\rho^{\Nbf}(r; \zbf)$ in both \texttt{GetRestartCandidate} and the restart condition are replaced with $\text{KKT}_{\omega}(\zbf)$.

\subsection{Computation of the normalized duality gap\label{app:normalized_duality_gap}}

The computation of the normalized duality gap follows the approach developed by \citet{applegate2023faster} for LPs and by \citet{xiong2024role} for conic programs. Below, we restate it in the context of our conic program formulation \eqref{eq:prob_cp}.

Note that  
\begin{equation}\label{pro:compute_rho}
\rho^{\Nbf}(r;\zbf):=\frac{1}{r}\cdot\left(\begin{array}{cc}
\max_{\hat{\zbf}} & \bbf^{\top}(\hat{\zbf}-\zbf)\\
\st & \hat{\zbf}\in\mcal X\times\mcal Y,\norm{\zbf-\hat{\zbf}}_{\Nbf}^{2}\leq r^{2}
\end{array}\right)
\end{equation}
where 
$\bbf=\left(\begin{array}{c}
\bbf_{1}\\
\bbf_{2}
\end{array}\right)=\left(\begin{array}{c}
\Gbf^{\top}\ybf-\cbf\\
\hbf-\Gbf \xbf
\end{array}\right).$  
For any $t$ satisfying $\|\zbf-\zbf(t)\|_{\Nbf}=r$, the optimizer of~\eqref{pro:compute_rho} is
\begin{equation}\label{pro:compute_zt}
\zbf(t) := \argmax_{\tilde{\zbf} \in \mcal X \times \mcal Y} \quad t \cdot \bbf^{\top}(\tilde{\zbf} - \zbf) - \frac{1}{2} \cdot \norm{\tilde{\zbf} - \zbf}_{\Nbf}^2.
\end{equation}
In Algorithm~\ref{alg:binary_search_normalized_gap}, Lines~\ref{alg:gap-doubling-start}--\ref{alg:gap-doubling-end} perform the doubling search, with $t_{k+1}=2t_k$ in Line~\ref{alg:gap-double}. Line~\ref{alg:gap-interval} sets the search interval $[t_{\text{left}},t_{\text{right}}]$ once $d_k>r$, with $\|\zbf-\zbf(t_{\text{left}})\|_{\Nbf}\leq r<\|\zbf-\zbf(t_{\text{right}})\|_{\Nbf}$. Lines~\ref{alg:gap-bisection-start}--\ref{alg:gap-bisection-end} then refine this interval by bisection. The argument of \citet[Lemma~A.2]{xiong2024role}, adjusted for the factor $\tfrac12$ in~\eqref{pro:compute_zt}, gives, for exact $\zbf(t)$ and $d_t:=\|\zbf(t)-\zbf\|_{\Nbf}\leq r$,
\[
0\leq\rho^{\Nbf}(r;\zbf)-\frac{\bbf^\top(\zbf(t)-\zbf)}{r}
\leq\frac{r^2-d_t^2}{2rt}\leq\frac{r}{2t}.
\]
Both searches return the gap estimate only when $d_t\leq r$ and this error bound is at most $\eps$. The doubling search can therefore terminate even if the distance never exceeds $r$.
\begin{algorithm}[htbp]
\caption{Bisection search for normalized duality gap $\rho^\Nbf(r;\zbf)$\label{alg:binary_search_normalized_gap}}
\small
\begin{algorithmic}[1]
    \Require{Initial bound $t_0>0$, iterate $\zbf\in\mcal X\times\mcal Y$, radius $r>0$, tolerance $\eps>0$}
    \State{$t_{\text{left}}\gets0$, $t_{\text{right}}\gets t_0$}
    \For{$k=0,1,\ldots$}\label{alg:gap-doubling-start} \Comment{Find a search interval}
        \State{$\tilde{\zbf}\gets\zbf(t_k)$, $d_k\gets\norm{\zbf-\tilde{\zbf}}_\Nbf$}
        \If{$d_k\leq r$ \textbf{and} $(r^2-d_k^2)/(2rt_k)\leq\eps$}
            \State{\Return $\bbf^\top(\tilde{\zbf}-\zbf)/r$}\label{alg:gap-bound-return} \Comment{Gap error bound satisfied}
        \EndIf
        \State{\textbf{if} $d_k>r$ \textbf{then set} $t_{\text{right}}\gets t_k$, set $t_{\text{left}}\gets t_k/2$ \textbf{if} $k>0$, \textbf{and break}}\label{alg:gap-interval}
        \State{$t_{k+1}\gets2t_k$}\label{alg:gap-double}
    \EndFor\label{alg:gap-doubling-end}
    \While{True}\label{alg:gap-bisection-start} \Comment{Refine by bisection}
        \State{$t_{\text{mid}}\gets(t_{\text{left}}+t_{\text{right}})/2$, $\tilde{\zbf}\gets\zbf(t_{\text{mid}})$}
        \State{$d_{\text{mid}}\gets\norm{\zbf-\tilde{\zbf}}_\Nbf$}
        \If{$d_{\text{mid}}\leq r$ \textbf{and} $(r^2-d_{\text{mid}}^2)/(2rt_{\text{mid}})\leq\eps$}
            \State{\Return $\bbf^\top(\tilde{\zbf}-\zbf)/r$} \Comment{Gap error bound satisfied}
        \EndIf
        \State{\textbf{if} $d_{\text{mid}}\leq r$ \textbf{then} $t_{\text{left}}\gets t_{\text{mid}}$ \textbf{else} $t_{\text{right}}\gets t_{\text{mid}}$}
    \EndWhile\label{alg:gap-bisection-end}
\end{algorithmic}
\end{algorithm}

Problem~\eqref{pro:compute_zt} reduces to two projections onto $\mathcal X$ and $\mathcal Y$ under the diagonal metric induced by $\Nbf$:
\begin{equation}\label{pro:compute_zt_2}
\zbf(t) = (\proj_{\mcal X}(\xbf+{t\tau}\cdot \bbf_{1}),\proj_{\mcal Y}\{\ybf+{t\sigma}\cdot \bbf_{2}\}) \ , 
\end{equation} 

\subsection{\label{sec:weight_adjustment}Primal weight adjustment}

As shown in Function~\ref{alg:adaptive_pdhg}, the primal weight is used to balance the primal and dual step sizes. A larger primal weight leads to a smaller primal step size $\tau$ and a larger dual step size $\sigma$. The update strategy follows a similar approach to that used in PDLP~\citep{applegate2021practical}, relying on the relative changes in the primal and dual variables.
Recent work by \citet{xiong2024accessible} provides partial theoretical support for primal-weight tuning, showing that a proper balance between the primal and dual iterates can accelerate restarted PDHG for LP problems.
In PDCS, the primal weight is updated only at each restart. The primal weight is initialized as follows:
\[
\texttt{InitializePrimalWeight}(\cbf,\hbf):=\begin{cases}
\frac{\|\cbf\|_2}{\|\hbf\|_2}, & \|\cbf\|_2,\|\hbf\|_2>10^{-10}\\
1 & \text{otherwise}
\end{cases}.
\]
The update at the beginning of a new restart is defined as:
\begin{equation*}
    \begin{aligned}
    &\texttt{PrimalWeightUpdate}(\zbf^{t,0},\zbf^{t-1,0},\omega^{t-1}):=\\
    &\begin{cases}
\exp\Brbra{\theta\log\brbra{\frac{\Delta_{y}^{t}}{\Delta_{x}^{t}}}+\brbra{1-\theta}\log\brbra{\omega^{t-1}}} & \Delta_{x}^{t},\Delta_{y}^{t}>10^{-10}\\
\omega^{t-1} & \text{otherwise}\ ,
\end{cases}
\end{aligned}
\end{equation*} 
where $\Delta_{x}^{t}=\norm{\xbf^{t,0}-\xbf^{t-1,0}}$ and $\Delta_{y}^{t}=\norm{\ybf^{t,0}-\ybf^{t-1,0}}$ denote the changes in the primal and dual variables, respectively. This update is designed to smooth changes in $\omega^t$ through the parameter $\theta$.
In practice, we set $\theta = 0.5$.
Occasionally, we observe that the solver can become unstable when the primal weight grows too large or too small, causing imbalanced updates between the primal and dual iterates. To mitigate this, if $\omega > 10^5$ or $\omega < 10^{-5}$, we reset $\omega$ to its initial value computed by \texttt{InitializePrimalWeight}.

\subsection{\label{subsec:Projection-onto-convex}Diagonal rescaling and projections onto rescaled cones}

Like many other first-order methods, the practical performance of PDCS is closely tied to the conditioning of the problem instance. One common strategy is diagonal rescaling, which may improve the condition number of the problem while preserving the sparsity pattern. This technique has become a standard step in many first-order solvers for LPs, including PDLP~\citep{applegate2021practical}, cuPDLP~\citep{lu2023cupdlp}, and ABIP~\citep{deng2024enhanced}. 

We rescale the constraint matrix $\Gbf \in \mbb{R}^{m \times n}$ to $\hat{\Gbf} = \Dbf_1 \Gbf \Dbf_2$, where $\Dbf_1$ and $\Dbf_2$ are positive diagonal matrices. The corresponding rescaled optimization problem is:
\begin{equation}\label{eq:rescaled_prob}
\min\ \   \inner{\hat{\cbf}}{\tilde{\xbf}} \ , \quad  
\st\ \   \hat{\mathbf{G}}\tilde{\xbf}-\hat{\mathbf{h}}\in\hat{\mathcal{K}}_{\textrm{d}}^{*} \ , \ \  \hat{\lbf}\leq(\tilde{\xbf})_{[n_{1}]}\leq\hat{\ubf} \ , \ \ (\tilde{\xbf})_{[n]\backslash[n_{1}]}\in\hat{\mathcal {K}}_{\textrm{p}} 
\end{equation}
where $\hat{\cbf}:=\Dbf_2\cbf$, $\hat{\Gbf}:=\Dbf_1\Gbf\Dbf_2$, $\hat{\hbf}:=\Dbf_1\hbf$, $\hat{\lbf}:=\Dbf_2^{-1}\lbf$, $\hat{\ubf}:=\Dbf_2^{-1}\ubf$, $\hat{\mcal K}_{\textrm{d}}^*:=\Dbf_1 \mcal K_{\textrm{d}}^*$, and $\hat{\mcal K}_{\textrm{p}}:=\Dbf_2^{-1}\mcal K_{\textrm{p}}$. 

The full PDCS framework with diagonal rescaling is summarized in Algorithm~\ref{alg:PDCS_after_precondition}.
\begin{algorithm}[h]
\caption{Primal-Dual Conic Programming Solver (PDCS) with preconditioning\label{alg:PDCS_after_precondition}}
    \begin{algorithmic}[1]
    \State{Compute rescaled instance data $\hat{\cbf},\hat{\Gbf},\hat{\hbf},\hat{\lbf}$ and $\hat{\ubf}$.}
    \State{Use Algorithm~\ref{alg:PDCS} to solve the rescaled problem~\eqref{eq:rescaled_prob} and obtain solution $\tilde{\xbf}$ and $\tilde{\ybf}$.}
    \State{Output solution for \eqref{eq:prob_cp}: $\xbf=\Dbf_2\tilde{\xbf}$ and $\ybf=\Dbf_1 \tilde{\ybf}$.}
    \end{algorithmic}
\end{algorithm}

This rescaling reformulates the problem in terms of non-standard cones $\hat{\mcal K}_{\textrm{p}}$ and $\hat{\mcal K}_{\textrm{d}}^*$. For linear programs (LPs), the rescaled nonnegative cones remain nonnegative, so projections remain straightforward. For rescaled second-order and exponential cones, the nontrivial projection cases reduce to scalar root-finding problems~\citep{lin2025practical}, solved by warm-started Newton iterations with bisection as a safeguard.

\section{Installation and use of PDCS}\label{sec:usage}
The source code of PDCS is publicly available at \url{https://github.com/ZikaiXiong/PDCS}. Currently, the solver supports the following types of cones:
\begin{itemize}
    \item Zero cone: $\zerobf^d$
    \item Nonnegative orthant: $\mbb R_+^d$
    \item Second-order cone: $\mcal K_{\soc}^{d+1}$
    \item Exponential cone: $\mcal K_{\exp}$
    \item Dual exponential cone: $\mcal K_{\exp}^*$
\end{itemize}
A rotated second-order cone is the image of a second-order cone under an invertible linear map, and models with rotated second-order cones are solved through this reformulation.

\subsection{Installation}
First clone the repository and install the Julia package:
\begin{codeblock}{bash}
    git clone https://github.com/ZikaiXiong/PDCS.git
    cd PDCS
\end{codeblock}
\begin{codeblock}{julia}
    using Pkg
    Pkg.develop(path="/path/to/PDCS")
\end{codeblock}

The CPU solver requires no CUDA installation. For the GPU solver, add \texttt{CUDA.jl} and compile the native projection kernels for the target GPU:
\begin{codeblock}{julia}
    using Pkg
    Pkg.add("CUDA")
\end{codeblock}
\begin{codeblock}{bash}
    make -C src/pdcs_gpu/cuda rebuild-gpu \
        CUDA_HOME=/path/to/cuda ARCH=sm_90
\end{codeblock}
The default architecture is \texttt{sm\_90}, which is used by NVIDIA H100 GPUs. The \texttt{ARCH} option should be changed for other GPU architectures. Detailed build and diagnostic instructions are provided in \texttt{src/pdcs\_gpu/cuda/README.md}.

\subsection{Usage}
PDCS provides JuMP~\citep{LubinDunningIJOC} and CVXPY~\citep{diamond2016cvxpy} interfaces, together with a direct Julia API. Importing the CPU solver does not load CUDA. The GPU solver is loaded only after \texttt{CUDA.jl} is imported.

\subsubsection{JuMP Usage\label{subsec:JuMP_Interface}}
The following example solves a small second-order cone program with the GPU solver:
\begin{codeblock}{julia}
    using CUDA
    using PDCS: PDCS_GPU
    using JuMP
    import MathOptInterface as MOI

    model = Model(PDCS_GPU.Optimizer)
    set_optimizer_attribute(model, "rel_tol", 1e-6)
    set_optimizer_attribute(model, "abs_tol", 1e-6)
    set_optimizer_attribute(model, "sparse_index_type", :auto)

    @variable(model, x[1:2])
    @variable(model, t)
    @objective(model, Min, t)
    @constraint(model, x[1] + x[2] == 1)
    @constraint(model, [t, x[1], x[2]] in SecondOrderCone())
    optimize!(model)

    status = termination_status(model)
    objective = objective_value(model)
    x_value = value.(x)

    metrics = MOI.get(
        JuMP.backend(model),
        MOI.RawOptimizerAttribute("result_metrics"),
    )
    solve_time_sec = metrics.solve_time_sec
    iterations = metrics.iterations
\end{codeblock} 
For the CPU implementation, use \texttt{using PDCS: PDCS\_CPU} and construct the model with \texttt{PDCS\_CPU.Optimizer}. With \texttt{sparse\_index\_type = :auto}, the GPU interface uses 32-bit sparse indices when the problem dimensions permit and 64-bit indices otherwise.

\subsubsection{Direct Julia API \label{subsec:Call_PDCS_from_Julia_directly}}
Users who already have a problem in standard conic form can call the solver directly and avoid JuMP/MOI model construction. The direct interface solves
\[
    \min_{\xbf}\ \cbf^\top\xbf
    \quad\text{s.t.}\quad \Gbf\xbf-\hbf\in\mcal K_{\textrm{d}}^{*},\qquad \lbf\leq\xbf\leq\ubf,
\]
where the rows of \(\Gbf\xbf-\hbf\) are ordered by zero, nonnegative, second-order, exponential, and dual exponential cones. For example:
\begin{codeblock}{julia}
    using CUDA
    using PDCS: PDCS_GPU

    solution = PDCS_GPU.rpdhg_gpu_solve(
        n = size(G, 2), m = size(G, 1), nb = length(c),
        c = c, G = G, h = h,
        mGzero = m_zero, mGnonnegative = m_nonnegative,
        socG = soc_dimensions,
        expG = number_exp, dual_expG = number_dual_exp,
        bl = lower_bounds, bu = upper_bounds,
        soc_x = Int[],
        rel_tol = 1e-6, abs_tol = 1e-6,
        time_limit = 3600.0,
    )

    x_value = Array(solution.x.recovered_primal.primal_sol.x)
    status = solution.info.exit_status
    solve_time_sec = solution.info.time
    iterations = solution.info.iter
\end{codeblock}
Here, \texttt{socG} lists the SOC dimensions, while \texttt{expG} and \texttt{dual\_expG} give the numbers of exponential and dual exponential cones. In this example, all structured cones are represented in \(\Gbf\xbf-\hbf\).
 
\subsubsection{CVXPY Interface \label{subsec:CVXPY_Interface}}
The CVXPY interface currently uses the GPU backend:
\begin{codeblock}{Python}
        import cvxpy as cp
        ...
        prob = cp.Problem(objective, constraints)
        prob.solve(solver=cp.PDCS, abs_tol=1e-6,
                   rel_tol=1e-6, verbose=True)
        ...
\end{codeblock}

\subsubsection{Termination Status}
Successful termination follows~\eqref{eq:stopping}. Table~\ref{tab:exit_code_and_status} lists the native PDCS exit codes and statuses.
\begin{table}[htbp]
    \centering
    \caption{Meanings of the exit codes and exit statuses \label{tab:exit_code_and_status}}
    \begin{tabular}{lcl}
        \toprule
        \textbf{Exit Status} & \textbf{Exit Code} & \textbf{Description} \\
        \midrule
        \texttt{:optimal} & 0 & Optimal solution found \\
        \texttt{:max\_iter} & 1 & Maximum number of iterations reached \\
        \texttt{:primal\_infeasible\_low\_acc} & 2 & Primal infeasible (low accuracy) \\
        \texttt{:primal\_infeasible\_high\_acc} & 3 & Primal infeasible (high accuracy)\\
        \texttt{:dual\_infeasible\_low\_acc} & 4 & Dual infeasible (low accuracy) \\
        \texttt{:dual\_infeasible\_high\_acc} & 5 & Dual infeasible (high accuracy) \\
        \texttt{:time\_limit} & 6 & Time limit reached \\
        \texttt{:continue} & 7 & Iteration continues (not terminated) \\
        \texttt{:numerical\_error\_nan} & 8 & NaN detected \\
        \texttt{:numerical\_error\_inf} & 9 & Infinite value detected \\
        \bottomrule
    \end{tabular}
\end{table}

\subsubsection{Default Settings}
By default, the GPU solver applies diagonal preprocessing, reflection, weighted averaging, adaptive restart, adaptive step sizes, and primal weight updates. The main GPU defaults are summarized in Table~\ref{tab:pdcs_default_parameters}. Restart and termination criteria are checked every 2,000 iterations by default.

\begin{table}[htbp]
    \centering
    \caption{Main default parameters of the GPU solver.}
    \label{tab:pdcs_default_parameters}
    \small
    \begin{tabular}{p{0.29\textwidth}p{0.15\textwidth}p{0.45\textwidth}}
        \toprule
        \textbf{Parameter} & \textbf{Default} & \textbf{Meaning} \\
        \midrule
        \texttt{rel\_tol}, \texttt{abs\_tol}
        & \texttt{1e-6}
        & Relative and absolute tolerances for the primal residual, dual residual, and objective gap. \\
        \texttt{time\_limit\_secs}
        & \texttt{1000.0}
        & Maximum solution time in seconds. \\
        \texttt{use\_scaling},\newline \texttt{rescaling\_method}
        & \texttt{true},\newline \texttt{:ruiz\_pock\_}\newline \texttt{chambolle}
        & Apply diagonal preprocessing using Ruiz and Pock--Chambolle scaling. \\
        \texttt{use\_adaptive\_restart}
        & \texttt{true}
        & Restart the iteration when the restart criterion is satisfied. \\
        \texttt{use\_adaptive\_step}
        & \texttt{true}
        & Adapt the step-size scale $\eta$ using Function~\ref{alg:adaptive_pdhg}. \\
        \texttt{use\_adaptive\_}\newline \texttt{step\_size\_weight}
        & \texttt{true}
        & Update the primal weight $\omega$ at each restart with smoothing parameter $\theta=0.5$. \\
        \texttt{max\_outer\_iter},\newline \texttt{max\_inner\_iter}
        & \texttt{10000}, \texttt{500000}
        & Set the outer- and inner-iteration limits. \\
        \texttt{sparse\_index\_type}
        & \texttt{:auto}
        & Select 32-bit or 64-bit GPU sparse indices from the problem size. \\
        \bottomrule
    \end{tabular}
\end{table}
In JuMP, parameters are changed through \texttt{set\_optimizer\_attribute}. In CVXPY, supported solver options are passed directly to \texttt{prob.solve}, as illustrated above for \texttt{abs\_tol} and \texttt{rel\_tol}. In the direct API, they are passed as keyword arguments, and \texttt{time\_limit} is unlimited by default.

\subsubsection{Projection Strategy}
Only SOCs use all four projection mappings. Algorithm~\ref{alg:projection-strategy} selects the SOC mapping from the number and dimensions of the cones, separately for the primal and dual SOC collections. Exponential, zero, and nonnegative cones use Thread-wise projection, with one thread per exponential cone and elementwise operations for zero and nonnegative cones. These rules also apply to mixed-cone models, with each cone type handled separately.

\begin{algorithm}[H]
\caption{Selection of the projection mapping for a collection of second-order cones}
\label{alg:projection-strategy}
\footnotesize
\begin{algorithmic}[1]
\Require SOC dimensions \(d_1,\ldots,d_l\), where \(l\geq1\) is the number of cones
\State \(d_{\max}\gets\max_i d_i\)
\State Select the base mapping \(B\) using
\Statex \hspace{\algorithmicindent}\(\displaystyle
B\gets
\begin{cases}
\mathrm{Thread}, & l=1,\ d_{\max}\leq3,\\
\mathrm{Warp},   & l=1,\ 4\leq d_{\max}\leq64,\\
\mathrm{Block},  & l=1,\ 65\leq d_{\max}<32{,}768,\\
\mathrm{Grid},   & l=1,\ d_{\max}\geq32{,}768,\\
\mathrm{Thread}, & l>1,\ d_{\max}\leq4,\\
\mathrm{Thread}, & l>1,\ 5\leq d_{\max}\leq64,\ l\geq1{,}024d_{\max},\\
\mathrm{Warp},   & l>1,\ 5\leq d_{\max}\leq64,\ l<1{,}024d_{\max},\\
\mathrm{Block},  & l>1,\ d_{\max}>64,\ (d_{\max}\geq1{,}024\ \lor\ l\leq998),\\
\mathrm{Warp},   & l>1,\ 64<d_{\max}<1{,}024,\ l>998.
\end{cases}\)
\If{\(B\in\{\mathrm{Block},\mathrm{Warp}\}\)}
    \State \(\mathcal T\gets\{i:d_i\leq4\}\), and \(\mathcal T\gets\varnothing\) if \(|\mathcal T|<256\)
    \State \(\mathcal W\gets\{i:5\leq d_i\leq32\}\), and \(\mathcal W\gets\varnothing\) if \(|\mathcal W|<64\)
    \State \(\mathcal R\gets\{1,\ldots,l\}\setminus(\mathcal T\cup\mathcal W)\)
    \State \(C\gets|\mathcal R|+\lceil|\mathcal T|/256\rceil+\lceil|\mathcal W|/8\rceil\)
    \If{\(\mathcal T\cup\mathcal W\neq\varnothing\) and \(5C\leq4l\)}
        \State Execute Thread on \(\mathcal T\), Warp on \(\mathcal W\), and \(B\) on \(\mathcal R\)
        \State \Return
    \EndIf
\EndIf
\State Execute \(B\) on the complete SOC collection
\end{algorithmic}
\end{algorithm}

The base rule balances parallelism within each cone against concurrency across cones. The second stage batches cones with similar dimensions to reduce load imbalance and scheduling overhead. These thresholds were calibrated on an NVIDIA H100 GPU, guided by the projection experiments in \citet{lin2025practical}, and may differ on other hardware. The environment variable \texttt{PDCS\_PROJECTION\_STRATEGY\_OVERRIDE} can be set to \texttt{thread}, \texttt{warp}, \texttt{block}, or \texttt{grid} to override the automatic SOC mapping.

\section{Some Experimental Results}\label{sec:experiments}

This section illustrates the performance of cuPDCS on three large-scale applications with different cone structures. Fisher market equilibrium involves many exponential cones, Lasso regression contains one large SOC, and multi-period portfolio optimization (MPO) has SOC constraints repeated across periods. All experiments were run on an NVIDIA H100 GPU with 80~GB of memory, with \texttt{abs\_tol} and \texttt{rel\_tol} set to $10^{-6}$ for the stopping test~\eqref{eq:stopping}. The results are from the experiments in \citet{lin2025practical}, which also compare cuPDCS with other solvers.

For each problem size, we generate five instances. Figure~\ref{fig:runtime_results} reports mean wall time, with error bars showing one standard deviation. 

\subsection{Fisher Market Equilibrium}\label{sec:fisher}

Let $m$ and $n$ be the numbers of buyers and goods. Buyer $i$ has budget $w_i$, $U_{ij}$ is the utility of good $j$ to buyer $i$, and $b_j$ is the supply of good $j$. Define the positive-valuation set $\mathcal E=\{(i,j):U_{ij}>0\}$. We solve the following exponential-cone formulation:
\begin{equation}\label{eq:fisher_conic}
\begin{aligned}
\min_{\{X_{ij}\}_{(i,j)\in\mathcal E},\pbf,\tbf}\quad
&-\sum_{i=1}^{m}w_i p_i\\
\text{s.t.}\quad
&\sum_{i:(i,j)\in\mathcal E}X_{ij}=b_j,\quad j=1,\ldots,n,\\
&t_i=\sum_{j:(i,j)\in\mathcal E}U_{ij}X_{ij},\quad (p_i,1,t_i)\in\mathcal K_{\exp},\quad i=1,\ldots,m,\\
&X_{ij}\geq0,\quad (i,j)\in\mathcal E.
\end{aligned}
\end{equation}
We generate utility matrices with density $0.2$, draw each nonzero entry of $\Ubf$ and each entry of $\wbf$ independently from $U[0,1]$, and set $b_j=0.25m$. Figure~\ref{fig:runtime_results}(a) reports the results for $n=1{,}000$.

cuPDCS solves all 45 reported instances, including problems with up to $800{,}000$ buyers.

\subsection{Lasso Regression}\label{sec:lasso}

For a data matrix $\Abf\in\mathbb R^{m\times n}$, response vector $\bbf\in\mathbb R^m$, and regularization parameter $\lambda>0$, the Lasso problem is
\begin{equation}\label{eq:lasso_original}
\min_{\xbf\in\mathbb R^n}\ \|\Abf\xbf-\bbf\|^2+\lambda\|\xbf\|_1.
\end{equation}
Writing $\xbf=\xbf_1-\xbf_2$ with $\xbf_1,\xbf_2\geq0$ gives the SOCP reformulation
\begin{equation}\label{eq:lasso_socp}
\begin{aligned}
\min_{r,\ybf,\xbf_1,\xbf_2}\quad
&2r+\lambda\onebf^\top(\xbf_1+\xbf_2)\\
\text{s.t.}\quad
&\ybf=\Abf(\xbf_1-\xbf_2)-\bbf,\qquad \xbf_1,\xbf_2\geq0,\\
&\left(\frac{1+r}{\sqrt{2}},\frac{1-r}{\sqrt{2}},\ybf\right)\in\mathcal K_{\mathrm{soc}}^{m+2}.
\end{aligned}
\end{equation}
Following~\citet{deng2024enhanced}, we generate $\Abf$ with sparsity $10^{-4}$ and draw its nonzero entries from $U[0,1]$. We set $\bbf=\Abf\widetilde{\xbf}+10^{-6}\onebf$, where $\widetilde{\xbf}$ is standard normal with half of its entries set to zero, and use $\lambda=\|\Abf^\top\bbf\|_\infty$. Figure~\ref{fig:runtime_results}(b) reports the results.

cuPDCS solves all 25 instances and scales to matrices with approximately $5.6\times10^8$ nonzero entries.

\subsection{Multi-period Portfolio Optimization}\label{sec:mpo}

Let $\wbf_{\tau+1}\in\mathbb R^{n+1}$ contain the weights of $n$ assets and cash at period $\tau+1$, and let $\wbf_b$ be a benchmark portfolio. The initial portfolio $\wbf_0$ is given. The subscript $[n]$ selects the first $n$ components, corresponding to the non-cash assets. Here, $\hat{\rbf}_{\tau+1}$ denotes predicted returns, $\hat{\mathbf\Sigma}_{\tau+1}\in\mathbb R^{n\times n}$ the estimated asset-return covariance matrix for period $\tau+1$, and $(\wbf_\tau^m)_{[n]}$ the market reference weights used in the market-neutrality constraint. We use the following constraint-based SOCP formulation:
\begin{equation}\label{eq:mpo_socp}
\begin{aligned}
\max_{\substack{\ubf_\tau\in\mathbb R^n,\,\wbf_{\tau+1}\in\mathbb R^{n+1}:\\ \tau=0,\ldots,T-1}}\quad
&\sum_{\tau=0}^{T-1}\hat{\rbf}_{\tau+1}^\top\wbf_{\tau+1}\\
\text{s.t.}\quad
&\onebf^\top(\wbf_{\tau+1}-\wbf_\tau)=0,\qquad \wbf_{\tau+1}\geq0,\\
&(\wbf_\tau^m)_{[n]}^\top\hat{\mathbf\Sigma}_{\tau+1}(\wbf_{\tau+1})_{[n]}=0,\\
&\ubf_{\tau,i}\geq\pm(\wbf_{\tau+1}-\wbf_b)_i,\quad i=1,\ldots,n,\\
&\sum_{i=1}^{n}\sqrt{(\hat{\mathbf\Sigma}_{\tau+1})_{ii}}\ubf_{\tau,i}\leq\gamma_{3\tau},\\
&\|\hat{\mathbf\Sigma}_{\tau+1}^{1/2}(\wbf_{\tau+1}-\wbf_b)_{[n]}\|\leq\gamma_{1\tau},\\
&-\gamma_{2\tau,i}\leq(\wbf_{\tau+1}-\wbf_\tau)_i\leq\gamma_{2\tau,i},\quad i=1,\ldots,n+1,
\end{aligned}
\quad \tau=0,\ldots,T-1.
\end{equation}
The objective maximizes cumulative predicted return. The constraints impose self-financing, long-only holdings, market neutrality, portfolio-risk limits, and transaction limits. The inter-period budget constraints neglect return-induced weight drift and transaction costs.

\begin{figure}[htbp]
  \centering
  \subfloat[Fisher market\label{fig:fisher_results}]{%
    \includegraphics[width=0.32\textwidth]{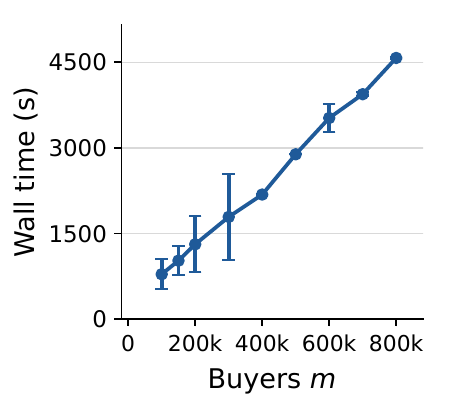}}\hfill
  \subfloat[Lasso\label{fig:lasso_results}]{%
    \includegraphics[width=0.32\textwidth]{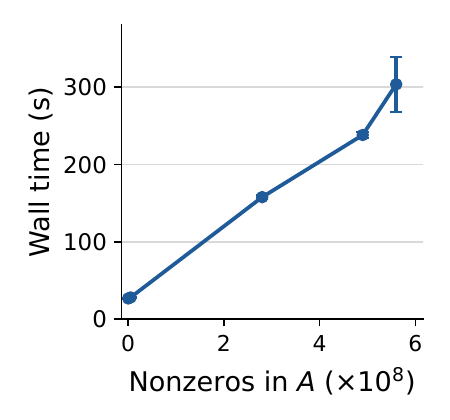}}\hfill
  \subfloat[MPO\label{fig:mpo_results}]{%
    \includegraphics[width=0.32\textwidth]{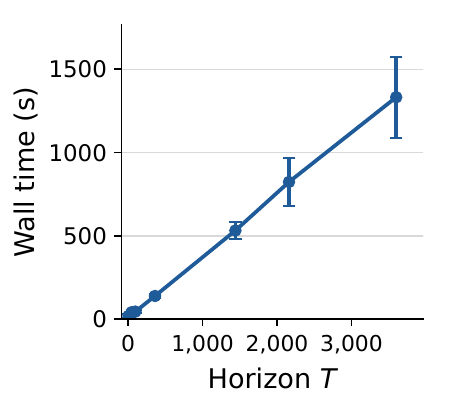}}
  \caption{cuPDCS runtime versus problem size with absolute and relative tolerances set to $10^{-6}$. Each point is the mean over five instances, and error bars show one standard deviation. All reported instances are solved.}\label{fig:runtime_results}
\end{figure}

We use 844 securities from the S\&P 500, NASDAQ-100, DOW 30, FTSE 100, and Nikkei 225. Data from 2024 determine the first-period returns, covariances, and prices, and subsequent periods are generated by adding controlled random noise. Five instances are generated for each horizon $T$. We set $\gamma_{1\tau}=\|\hat{\mathbf\Sigma}_{\tau+1}^{1/2}(\wbf_{1/n}-\wbf_b)_{[n]}\|$, $\gamma_{2\tau,i}=0.05$, and $\gamma_{3\tau}=0.05$, where $\wbf_{1/n}$ assigns weight $1/n$ to each asset and zero to cash. Figure~\ref{fig:runtime_results}(c) reports the results.

cuPDCS solves all 35 instances. The mean runtime increases with the horizon, reaching $1{,}332$ seconds at $T=3{,}600$.

\section{Concluding Remarks}
PDCS combines matrix-free primal-dual iterations with GPU cone projections for large-scale conic optimization. The three applications illustrate its ability to handle different cone structures at moderate accuracy. JuMP and CVXPY provide modeling interfaces, while the direct Julia API accepts conic data without an intermediate modeling layer. PDCS complements interior-point solvers, which remain a natural choice for smaller problems requiring high accuracy.

\renewcommand \thepart{}
\renewcommand \partname{}

\bibliographystyle{abbrvnat}
\bibliography{ref}

\end{document}